\documentclass[reqno]{amsart}

\usepackage{amsmath,amssymb,amsthm}

\usepackage{tikz}
\usepackage{caption}
\usepackage{graphicx}
\usepackage{graphics}

\usepackage{hyperref}

\newtheorem*{thm}{Theorem}

\usepackage{url}
\usepackage{xcolor}

\begin{document}

\title[]{A Remark on the Odd Area of Unit Disks }

\author[]{Stefan Steinerberger}

\address{Department of Mathematics and Department of Applied Mathematics, University of Washington, Seattle, WA 98195, USA}
 \email{steinerb@uw.edu}

\begin{abstract}  Let $\mathcal{F}$ be a family of $n$ unit disks in $\mathbb{R}^2$ with $n$ being odd.  We use $\mbox{OA}(\mathcal{F})$ to denote the area of the set of points that is covered by an odd number of disks.  
The purpose of this note is to disprove the conjecture $\mbox{OA}(\mathcal{F}) \geq \pi$ which was suggested in the literature and to present some examples.
\end{abstract}
\maketitle

The following problem has been discussed a number of times \cite{car, ore, pin0, pin}: given a collection $\mathcal{F}$ of $n$ unit disks in $\mathbb{R}^2$, is it true that the measure  $\mbox{OA}(\mathcal{F})$ of the set of points that is covered by an odd number of disks is at least $\pi$? Taking all the disks to be identical shows that this would be optimal (in particular, it is stated as Problem A in \cite{car} and described as \textit{the most interesting problem related to odd area} in \cite{ore}).
Understanding the best constant in $\mbox{OA}(\mathcal{F}) \geq c$ is also an open problem in Pak's textbook \cite[Exercise 15.14]{pak}. It is not even clear whether $c>0$.

\begin{thm} Let $\mathcal{F}$ be six disks whose centers are equispaced on the boundary of the equilateral triangle with length 1/2 (three in the vertices, three in the midpoints between vertices)  with a seventh disk placed in the center of the triangle. Then 
$$\emph{OA}(\mathcal{F}) \leq 3.12.$$
\end{thm}
\vspace{-10pt}

\begin{center}
\begin{figure}[h!]
\begin{tikzpicture}
\node at (0,0) {\includegraphics[width=0.3\textwidth]{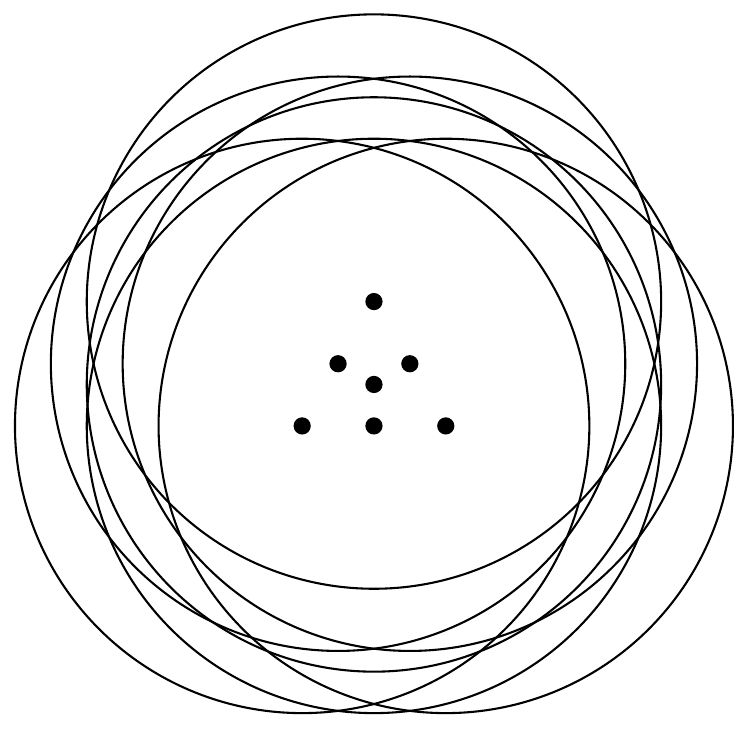}};
\node at (6,0) {\includegraphics[width=0.3\textwidth]{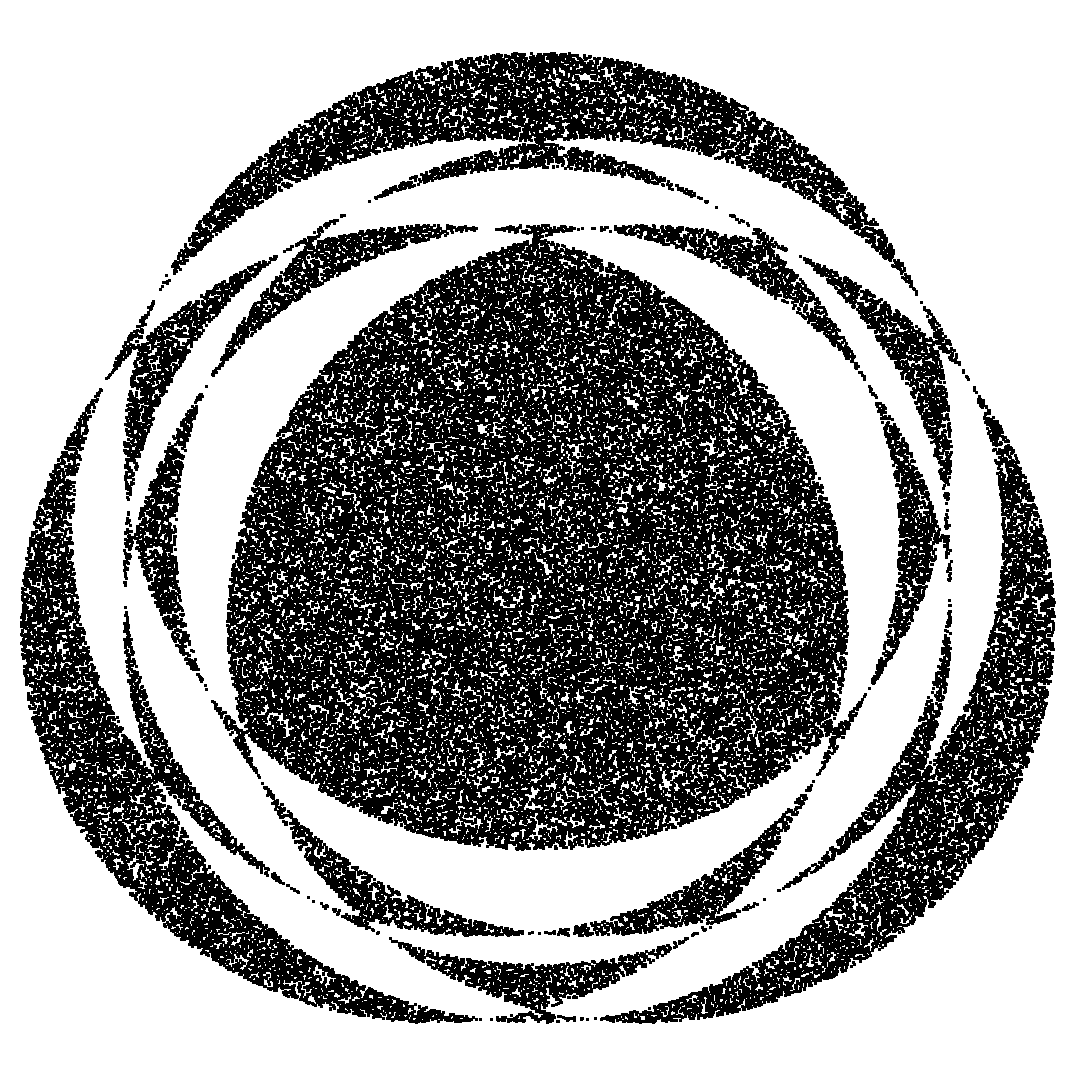}};
\end{tikzpicture}
\caption{Left: the centers of the 7 disks. Right: the points in the plane that are contained in an odd number of disks (by sampling random points and checking whether they are in the set).}
\end{figure}
\end{center}
\vspace{-10pt}

This also disproves Conjecture 1 in \cite{pin}. As for the proof, the area in question could presumably be computed in closed form. A more mundane numerical grid counting method gives area $\sim 3.108$ (also suggested by Monte Carlo methods); in light of subsequent examples, the precise value appears less important since $\mbox{OA}(\mathcal{F}) \leq 2.5$ seems possible. We discuss two cases, $n=51$ and $n=151$, which were investigated using Google DeepMind's \textsc{AlphaEvolve} \cite{novikov}. \newpage

\textbf{51 disks.} It appears that for a `small' number of disks (say $n \leq 100$), there exist families that are roughly `round' in shape in the sense that the centers of the disks appear to be distributed in a disk-like region. Using grid-based methods to approximate the area appears to be impossible since an extremely fine mesh would be required.  The Monte Carlo method is oblivious to the complexity of the shape of the domain and suggests, for the example shown in Figure 2, that
$$\mbox{OA}(\mathcal{F}) \leq 2.73$$

\begin{center}
\begin{figure}[h!]
\begin{tikzpicture}
\node at (0,0) {\includegraphics[width=0.3\textwidth]{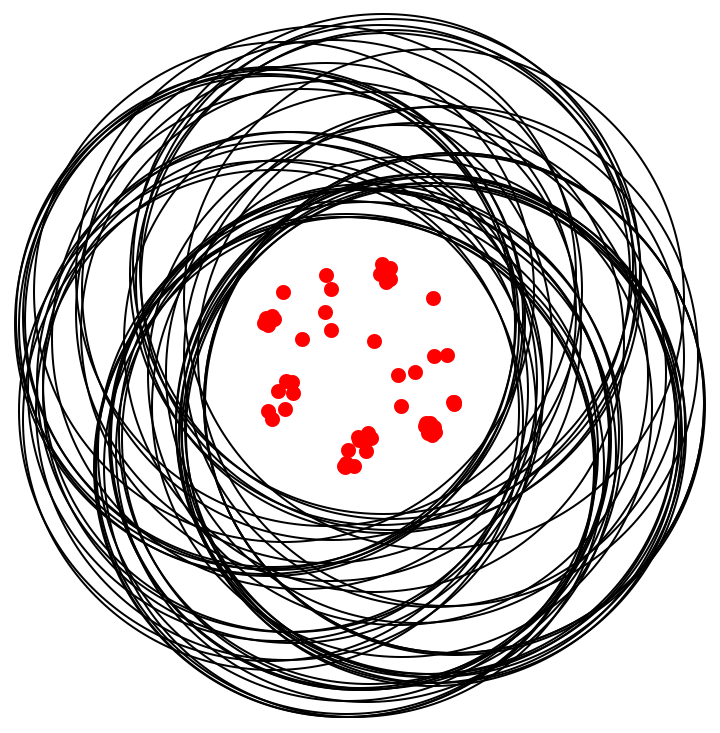}};
\node at (4,0) {\includegraphics[width=0.3\textwidth]{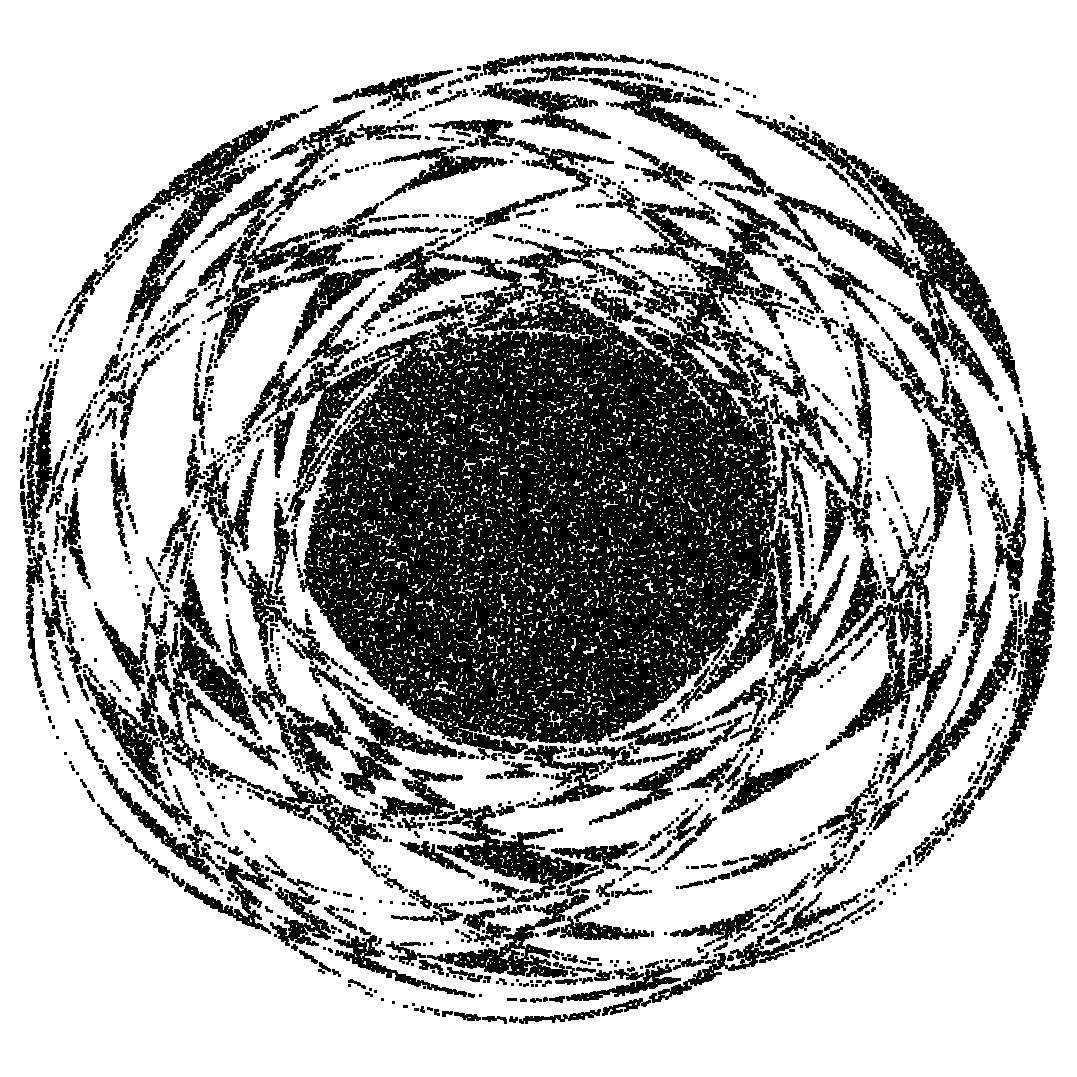}};
\node at (8.2,0) {\includegraphics[width=0.24\textwidth]{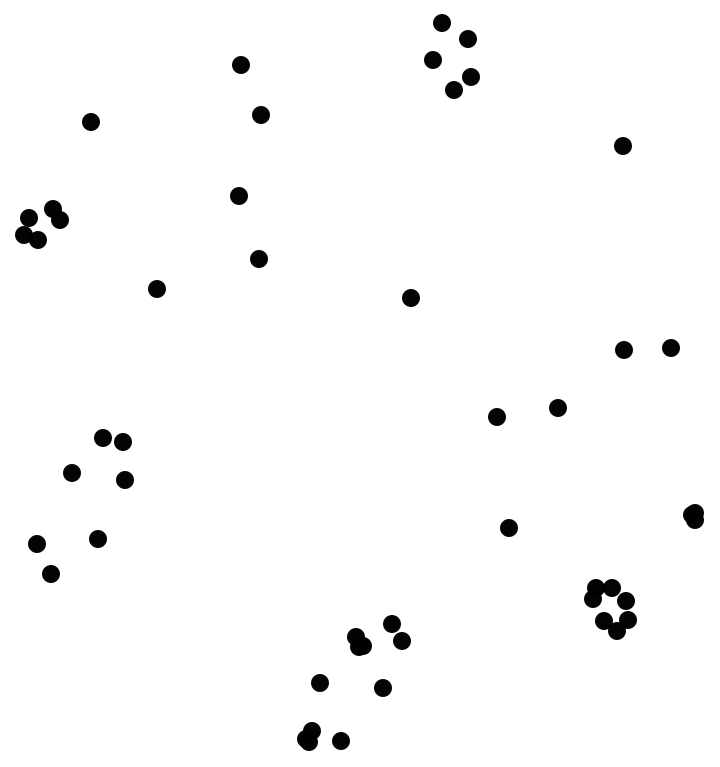}};
\end{tikzpicture}
\caption{Left: centers of the 51 disks (in red) and the corresponding unit circles. Middle: the region that is contained in an odd number of disks (by sampling random points and checking whether they are in the set). Right: the centers of the circles.}
\end{figure}
\end{center}
 \textbf{151 disks.} As the number of disks increases there is a change in shape of the local minimizers identified by \textsc{AlphaEvolve}: the centers of the disk appear to lie in an elongated ellipse-like convex domain and on a type of curve going through it. The example shown in Fig. 3 appears to satisfy, again via Monte-Carlo, that 
 $$\mbox{OA}(\mathcal{F}) \leq 2.51.$$ 
 It is completely unclear whether this particular configuration is in any way indicative of a larger family of examples, whether it is specific to $n=151$, or whether it is an artifact of the optimization method.

 \begin{center}
\begin{figure}[h!]
\begin{tikzpicture}
\node at (0,0) {\includegraphics[width=0.35\textwidth]{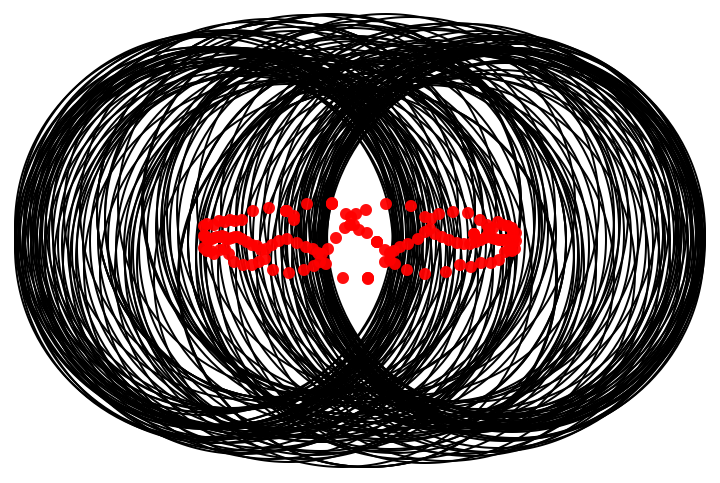}};
\node at (6.5,0) {\includegraphics[width=0.55\textwidth]{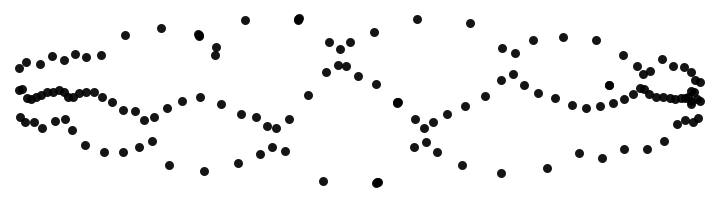}};
\end{tikzpicture}
\caption{left: centers of the 151 disks (in red) and the corresponding unit circles. Right: zoom into the centers of the circles.}
\end{figure}
\end{center}

\textbf{More disks.} We were unable to produce results of a comparable quality with more disks.  We were able to find examples of 201 disks for which $\mbox{OA}(\mathcal{F}) \leq 3.07$ but these arise as a perturbation of the case where all the disks are in the same place and do not seem to provide any additional insight; since disks can occupy the same space, we know that there exist configurations with $\mbox{OA}(\mathcal{F}) \leq 2.51$ for all $n \geq 151$, however, we did not find any such examples numerically (presumably because the dramatically increasing combinatorial complexity which makes optimization more complicated).
Indeed, it appears that the optimization landscape is \textit{incredibly} subtle. Reasons may include
\begin{enumerate}
\item $n$ disks force one to consider $2^n$ potential intersection patterns.
\item Monte-Carlo sampling is (presumably) the only viable approach which makes standard gradient-based techniques with a small stepsize more complicated since there is a random error on top of every evaluation.
\item The functional $\mbox{OA}$ as a function of the $2n$ coordinates is very unstable; many of the intersection areas are thin, moving points even very slightly can change the odd area a lot. The example for $n=151$ shown above, resulting in  $\mbox{OA}(\mathcal{F}) \sim 2.51$ degenerates to  $\mbox{OA}(\mathcal{F}_*) \sim 2.72$, where $\mathcal{F}_*$ is $\mathcal{F}$ with coordinates rounded to three digits.
\item Finally, the configuration where all disks are in the same space appears to have a large basin of attraction.
\end{enumerate}
While the inequality $\mbox{OA}(\mathcal{F}) \geq \pi$ would have been nice, its failure seems to happen in a \textit{most} intriguing way, both in terms of the complexity of the counterexamples as well as the optimization problem (which may be a useful benchmark problem in its own right). Could it be that $\mbox{OA}(\mathcal{F}) \leq \varepsilon$ is possible for every $\varepsilon >0$?\\

  \textbf{Statement on AI use.} The Theorem is due to the author (inspired by numerics for $n=7$).
  The optimization was carried out using Google DeepMind's \textsc{AlphaEvolve} \cite{novikov} which was prompted with a problem
  description and a large number of additional constraints to avoid degeneracies (centers of all the disks have to lie in $[-3,3]^2$, no two centers are allowed to be too close, at least two centers are at least a certain distance from each other); constraints appear to be necessary to avoid the centers all collapsing into a single point (and, once $n \geq 151$, we have been unable to find constraints that would prevent that from happening -- we take this to be an indication of the growing combinatorial complexity of the problem rather than any fundamental change in the overall behavior).  The results were independently verified by the author using \textsc{Mathematica} which was also used to create the figures. \\

 \textit{Coordinates of the $n=51$ example. } 
$$(-0.232, 0.155), (0.056, 0.145), (0.123, 0.394), (0.046, -0.243), (0.289, -0.232)$$
$$(-0.349, 0.245), (-0.376, 0.236), (-0.367, -0.133), (0.024, -0.297), (0.12, 0.438)$$
$$(-0.024, -0.356), (-0.298, -0.128), (-0.136, 0.409), (0.274, -0.221)$$
$$(0.001, -0.248), (0.166, -0.115), (-0.063, -0.354), (0.091, 0.456), (0.377, -0.098)$$
$$ (-0.114, 0.352), (0.034, -0.223), (0.104, 0.38), (0.377, -0.106), (0.296, 0.317)$$
$$ (-0.328, -0.053), (-0.006, -0.239), (-0.056, -0.344), (-0.116, 0.189)$$
$$(-0.351, -0.167), (0.08, 0.414), (-0.382, 0.216), (-0.268, -0.061), (0.222, 0.021)$$
$$ (0.296, 0.086), (-0.366, 0.21), (0.299, -0.197), (0.265, -0.183), (-0.047, -0.29)$$
$$ (-0.139, 0.26), (-0.293, -0.014), (0.35, 0.089), (-0.341, 0.233), (-0.306, 0.344)$$
$$(0.283, -0.183), (0.374, -0.1), (0.153, 0.011), (-0.27, -0.017), (-0.003, -0.25)$$
 $$  (0.302, -0.219), (0.262, -0.196), (-0.059, -0.357)\\$$

\textit{Coordinates of the $n=151$ example. } As
mentioned above, truncating the coordinates to three digits leads to a degeneration
of the example from $\mbox{OA}(\mathcal{F}) \sim 2.5$ to $\mbox{OA}(\mathcal{F}) \sim 2.72$.
For the sake of brevity, we only list the truncated example; the reader can
find an attached notebook with the original example in high precision.
$$(-0.706, 0.001), (0.793, -0.018), (0.693, 0.062), (-0.192, -0.131), (0.562, 0.137)$$
$$ (0.023, 0.157), (0.121, -0.122), (-0.691, 0.008), (-0.493, 0.166), (0.127, 0.19)$$
$$ (-0.064, 0.078), (0.08, -0.015), (0.793, 0.013), (-0.36, 0.122), (0.244, -0.023)$$
$$ (-0.515, -0.106), (-0.101, -0.204), (0.462, -0.004), (0.443, -0.171), (0.722, 0.091)$$
$$ (-0.085, 0.133), (0.149, -0.11), (0.333, 0.118), (-0.769, 0.011), (-0.135, 0.005)$$
$$ (-0.034, 0.134), (-0.388, -0.179), (0.03, -0.206), (0.806, -0.004), (0.792, -0.01)$$
$$(0.408, 0.137), (0.177, -0.133), (0.801, 0.042), (0.367, 0.106), (0.028, 0.031)$$
$$ (-0.401, 0.149), (-0.261, -0.049), (-0.473, -0.164), (0.081, -0.013), (0.814, 0.035)$$
$$ (-0.299, -0.04), (0.66, 0.075), (-0.397, 0.), (-0.4, 0.151), (-0.346, -0.017)$$
$$ (0.777, -0.003), (-0.728, 0.012), (-0.673, 0.096), (0.594, 0.029), (-0.092, 0.059)$$
$$ (-0.719, 0.), (0.167, -0.06), (-0.786, 0.079), (-0.555, -0.033), (-0.707, -0.079)$$
$$ (-0.584, -0.134), (-0.579, 0.15), (0.122, -0.054), (0.294, 0.003), (0.758, -0.067)$$
$$(-0.533, -0.055), (0.651, 0.006), (0.767, -0.004), (-0.635, -0.001), (0.033, -0.207)$$
$$(0.029, -0.208), (0.679, 0.02), (-0.015, 0.05), (-0.739, 0.017), (-0.782, 0.005)$$
$$ (-0.809, -0.005), (-0.756, 0.099), (0.726, -0.107), (-0.547, -0.122), (-0.675, 0.012)$$
$$ (-0.585, -0.032), (-0.8, -0.062), (0.503, -0.02), (-0.306, -0.161), (-0.158, 0.19)$$
$$ (-0.754, 0.011), (-0.728, 0.09), (0.675, 0.054), (0.079, -0.012), (-0.83, 0.019)$$
$$ (-0.631, -0.134), (0.796, -0.061), (-0.821, -0.06), (-0.782, -0.076), (-0.75, -0.059)$$
$$(-0.796, 0.001), (0.52, -0.137), (0.724, -0.001), (-0.837, 0.069), (-0.159, 0.186)$$
$$ (0.814, -0.01), (-0.819, 0.084), (-0.835, -0.048), (-0.656, 0.013), (-0.252, -0.139)$$
$$ (-0.701, 0.105), (-0.213, -0.075), (-0.362, 0.101), (-0.182, -0.053), (-0.045, 0.074)$$
$$ (0.755, -0.005), (0.63, -0.004), (0.81, -0.05), (0.69, 0.007), (-0.058, 0.116)$$
$$ (-0.441, -0.01), (-0.677, -0.117), (-0.161, 0.189), (0.785, -0.001), (0.74, -0.003)$$
$$ (0.257, 0.178), (0.539, -0.027), (0.631, -0.126), (-0.51, -0.049), (-0.29, 0.185)$$
$$(0.776, 0.071), (-0.402, 0.152), (-0.819, -0.003), (0.594, 0.03), (0.706, 0.)$$
$$ (0.792, 0.06), (0.332, -0.186), (0.422, 0.01), (0.685, -0.126), (0.748, 0.075)$$
$$ (0.603, -0.014), (-0.235, -0.071), (0.387, 0.028), (0.572, -0.022), (0.237, -0.164)$$
$$ (-0.478, -0.027), (0.799, 0.011), (0.577, -0.148), (0.201, -0.042), (0.777, -0.056)$$
$$ (-0.224, -0.121), (0.482, 0.146), (0.361, 0.054), (0.627, 0.102), (0.667, 0.021), (-0.638, 0.101)$$
$$ (0.146, -0.075), (-0.61, -0.012), (0.333, 0.041), (-0.837, 0.016), (-0.725, -0.054)$$

\end{document}